\documentclass{amsart}

\def\C{\mathbb{C}}
\def\Q{\mathbb{Q}}
\def\Z{\mathbb{Z}}
\def\Li{\operatorname{Li}}
\newcommand{\EIS}[1]{\emph{EIS}~#1}
\newtheorem{theorem}{Theorem}
\newtheorem{proposition}{Proposition}
\newcounter{noteno}
\newenvironment{note}{\refstepcounter{noteno}%
	\medbreak\noindent{\bf Note~\thenoteno.}}%
	{\hfill{\Large $\lhd$}\medbreak}

\begin{document}

\title[Non-holonomy]{\bf On the non-holonomic character of logarithms,
powers, and the $n$th
prime function}
\author{Philippe Flajolet \and Stefan Gerhold \and Bruno Salvy}
\address{Philippe Flajolet: Algorithms Project, INRIA Rocquencourt, F-78153 Le
Chesnay (France)\\
\protect{\tt Philippe.Flajolet AT inria.fr}}

\address{Stefan Gerhold: Research Institute for Symbolic Computation,
Johannes Kepler University Linz (Austria)
\\
\protect{\tt stefan.gerhold AT risc.uni-linz.ac.at}}
\address{Bruno Salvy: Algorithms Project, INRIA Rocquencourt, F-78153 Le
Chesnay (France)\\
\protect{\tt Bruno.Salvy AT inria.fr}}

\date{January 21, 2005}

\begin{abstract}
We establish that the sequences formed by logarithms and by ``fractional''
powers of integers, as
well as
the sequence of prime numbers, are non-holonomic, thereby answering three
open problems of Gerhold [\emph{El. J. Comb.} {\bf 11} (2004), R87].
Our proofs depend on basic complex analysis, namely
a conjunction of
the Structure Theorem for singularities of solutions to
linear differential equations and of an Abelian theorem.
A brief discussion is offered regarding the scope of
singularity-based methods and several naturally
occurring sequences are proved to be non-holonomic.
\end{abstract}

\maketitle

\begin{small}
\begin{flushright}
\sl Es ist eine Tatsache, da{\ss} die genauere Kenntnis\\ des Verhaltens
einer analytischen Funktion\\ in der N{\"a}he ihrer singul\"aren
Stellen\\
eine Quelle von arithmetischen S{\"a}tzen ist.\footnote{%
        \em ``It is a fact that the precise knowledge of the
        behaviour of an analytic function in the 
        vicinity of its singular points is a source of arithmetic properties.''
}\\[2mm]
{\sc --- Erich Hecke} \rm\cite[Kap.~VIII]{Hecke23}
\end{flushright}
\end{small}

\section*{Introduction}

A  sequence $(f_n)_{n\ge0}$
of complex numbers is said to be \emph{holonomic} (or
\emph{$P$-recursive})
if it satisfies a linear recurrence with coefficients that are
polynomial in the index~$n$, that is,
\begin{equation}\label{h1}
p_0(n)f_{n+d}
+p_1(n)f_{n+d-1}+\cdots+p_d(n)f_{n}=0, \qquad
n\ge 0,
\end{equation}
for some polynomials $p_j(X)\in\C[X]$. A formal power series 
$f(z)=\sum_{n\ge0}f_n z^n$ is \emph{holonomic} (or
\emph{$\partial$-finite}) if it satisfies a linear differential
equation
with coefficients that are polynomial in the variable~$z$, that is,
\begin{equation}\label{h2}
q_0(z)\frac{d^e}{dz^e}
f(z)+q_1(z)\frac{d^{e-1}}{dz^{e-1}} f(z)+\cdots +q_e(z) f(z)=0.
\end{equation}
for some polynomials $q_k(X)\in\C[X]$. It is well known that a sequence
is holonomic if and only if its generating series is
holonomic. (See Stanley's book~\cite{Stanley98} for the basic
properties of these sequences and series.)  
By extension, a function analytic at~0 is called holonomic if its
power series representation is itself holonomic. 

Holonomic sequences encapsulate many of the combinatorial sequences of common
interest, for instance, a wide class of sums involving binomial coefficients.
At the same time, they enjoy a varied set of closure
properties and 
several formal mechanisms have been recognised
to lead systematically 
to holonomic sequences---what we have in mind here includes
finite state models and regular
grammars leading to rational (hence, holonomic) functions, context-free
specifications leading to algebraic (hence, holonomic) functions, 
a wide class of problems endowed with symmetry conducive to holonomic
functions (via Gessel's theory~\cite{Gessel90}). On these aspects, we
may rely on general references like~\cite{FlSe04,Stanley98} as well as
on many works of Zeilberger, who is to be held accountable for 
unearthing the power of the holonomic framework accross
combinatorics; see~\cite{PeWiZe96,Zeilberger90}. 
Thus, in a way, a non-holonomicity result
represents some sort of a structural complexity lower bound.

This note answers three problems described as open
in an article of Stefan
Gerhold~\cite{Gerhold04} very recently published in the
\emph{Electronic Journal of Combinatorics}. 

\begin{proposition}\label{log-thm}
The sequence $f_n=\log n$ is not holonomic.
\end{proposition}
\noindent
(For definiteness, we agree that $\log 0 \equiv 0$ here.)

For~$\alpha$ an integer, the sequence
\[
h_n=n^\alpha\]
is clearly holonomic. (As a matter of fact, the generating function is
rational if $\alpha\in\Z_{\ge0}$ and of polylogarithmic type if
$\alpha\in\Z_{<0}$.)
Gerhold~\cite{Gerhold04} proved that for any $\alpha$ that is rational but not
integral, $h_n$ fails to be holonomic. For instance, $h_n=\sqrt{n}$
fails to be holonomic because, in essence,
$\Q(\sqrt{2},\sqrt{3},\ldots)$ is not a finite extension
of~$\Q$.

\begin{proposition}\label{pow-thm}
For $\alpha\in\C$, the sequence of powers $h_n=n^\alpha$ is holonomic
if and only if $\alpha\in\Z$.
\end{proposition}
\noindent
(We agree that $h_0=0$.)

\begin{proposition}\label{prime-thm}
The sequence $g_n$ defined by the fact that $g_n$ is the $n$th prime
is non-holonomic. 
\end{proposition}
\noindent
(We agree that $g_0=1$, $g_1=2$, $g_2=3$, $g_3=5$, and so on.)

\smallskip

Proposition~\ref{log-thm}, conjectured by Gerhold in~\cite{Gerhold04} 
was only proved
under the assumption that a difficult conjecture of number
theory (Schanuel's conjecture) holds. 
The author of~\cite{Gerhold04}
describes the statement of our
Proposition~\ref{pow-thm} as a \emph{``natural conjecture''}.
Proposition~\ref{prime-thm}
answers an explicit question of Gerhold who writes: \emph{``we do not know of
any proof that the sequence of primes is non-holonomic''.}  

\smallskip

Our proofs are plainly based on the combination of two facts.
First, holonomic objects satisfy rich
closure properties. In particular, we make use of closure of under sum, product and composition with an algebraic function.
Second, the asymptotic behaviour of
holonomic sequences, which is reflected by the asymptotic behaviour at
singularities of their generating functions, is rather strongly
constrained. For instance, iterated logarithms or negative powers of
logarithms are ``forbidden'' and estimates like
\begin{equation}\label{e0}
a_n \mathop{\sim}_{n\to+\infty} \log\log n, \qquad
b(z)\mathop{\sim}_{z\to1-} \frac{1}{\log(1-z)},
\end{equation}
are sufficient to conclude that 
the sequence $(a_n)$ and the function $b(z)$ are non-holonomic.
(See details below.)
A conjunction of the previous two ideas then 
perfectly describes the strategy of this note: 
In order to prove that a sequence is
non-holonomic, it suffices to locate a ``derived sequence'' (produced
by holonomicity-preserving transformations) that exhibits a suitable
combination of kosher asymptotic terms with a foreign
non-holonomic element, like in~(\ref{e0}).

We choose here to operate directly with generating
functions. Under this scenario, one can rely on the well established
classification of singularities elaborated at the end of the
nineteenth century by Fuchs~\cite{Fuchs66}, Fabry~\cite{Fabry85}
and others. A summary of what is known is 
found in standard
treatises, for instance the ones by Wasow~\cite{Wasow87} and
Henrici~\cite{Henrici77}. What we need of this theory is summarised by
Theorem~\ref{classif-thm} of the next section.

The relation between asymptotic behaviour of sequences and local
behaviour of the generating functions is provided by a classical
Abelian theorem, stated as 
Theorem~\ref{abel-thm} below.

\smallskip

\begin{note}
A methodological remark is in order at this stage. A glance
at~(\ref{e0}) suggests two possible paths for proving
non-holonomicity: one may \emph{a priori}
operate equally well with sequences
or with generating functions.  The latter is what we have opted to do here.
The former approach with sequences seems workable, but it requires a strong 
structure theorem analogous to Theorem~\ref{classif-thm} below
for recurrences, i.e., difference equations. 
An ambitious programme towards such a goal was
undertaken by Birkhoff and 
Trjitzinsky~\cite{Birkhoff38,BiTr30} in the 1930's, their works being later 
followed by Wimp and
Zeilberger in~\cite{WiZe85}. However, what is available 
in the classical literature is largely
a set of \emph{formal} solutions to difference equations and recurrences, 
and the relation of these to
actual (analytic) solutions represents a difficult problem 
evoked in~\cite[p.~168]{WiZe85} and~\cite[p.~1138]{Odlyzko95};
see also~\cite{BrFaIm00} for recent results relying on multisummability. 
\end{note}

\begin{note}
In this short paper, we do nothing but
assemble some rather well-known
facts of complex analysis, and logically organise them towards 
the goal of proving certain sequences to be non-holonomic. Our
purpose is thus essentially pedagogical.  
As it should become transparent soon, a rough heuristic in
this range of problem is the following: \emph{Almost
anything is
non-holonomic unless it is holonomic by design}. 
(This na{\"\i}ve remark cannot of course be universally true and there are
surprises, e.g., some sequences may eventually admit algebraic or holonomic
descriptions for rather deep reasons. Amongst such cases,
we count the enumeration
of~$k$-regular graphs and various types of
maps~\cite{Gessel90,GoJa83},
the enumeration of permutations with bounded-length increasing
subsequences, the Ap\'ery sequence~\cite{Poorten79}
related to a continued fraction
expansion of~$\zeta(3)$,
as well as the appearance of holonomic functions in the theory of modular
forms, for
which we refer to the beautiful exposition of Kontsevich-Zagier~\cite{KoZa01}.) 
\end{note}

\section{Methods}

From the most basic theorems regarding the existence of analytic
solutions to differential equations (e.g., \cite[Th.~9.1]{Henrici77}),
any function~$f(z)$ analytic at~0 
that is holonomic can be continued analytically along any path that avoids the
finite set~$\Sigma$ of points defined as roots of the equation $p_0(z)=0$,
where $p_0$ is the leading coefficient in~(\ref{h1}). 
Figuratively:

\begin{theorem}[Finiteness of singularities]\label{finite-thm}
A holonomic function has only finitely many singularities.
\end{theorem}

This theorem gives immediately as non-holonomic a number of 
sequences enumerating classical combinatorial structures.
\begin{itemize}
\item[---] Integer partitions, whose generating function
is $P(z)=\prod(1-z^n)^{-1}$, as the function admits the unit circle as a natural
boundary. The same argument applies to integer partitions with
summands restricted to any infinite set (e.g., primes), partitions
into distict summands, plane partitions, and so on.
More generally, combinatorial classes defined by an unlabelled set or
multiset construction~\cite{FlSe04} are non-holonomic, unless a rather drastic 
combinatorial simplification occurs.
\item[---] Alternating (also known as zig-zag, up-and-down, cf\footnote{%
	In order to keep this note finite, we refer to some of the
	combinatorial problems by means of their number in
	Sloane's \emph{Encyclopedia
	of Integer Sequences (EIS)}, see~\cite{Sloane00}.
}~\EIS{A000111})
permutations
with exponential generating function $\tan z + \sec z$, as they have
the odd multiples of $\frac\pi2$ as set of poles\footnote{%
	Stanley~\cite{Stanley80} describes an algebraic proof
	dependent on the fact that $\exp(z)$ is nonalgebraic (his
	Example~4.5), then goes on to observe in his~\S4.a
	that $\sec z$ ``\emph{has infinitely many poles}''.}. 
A similar argument
applies to preferential arrangements (also known as ordered set
partitions or surjections, cf \EIS{A000670}), Bernoulli numbers, and the like.
\item[---] Necklaces (equivalently Lyndon words,
irreducible polynomials), whose generating function admits the unit
circle as a natural boundary. More generally, ``most'' unlabelled
cycles are non-holonomic.
\item[---] Unlabelled plane trees (\EIS{A000081}), whose implicit specification involves an
unlabelled multiset construction. 
\end{itemize}

In many cases, the criterion above is too brutal. For instance it does
not preclude holonomicity for the Cayley tree function,
\begin{equation}\label{cayley}
T(z)=\sum_{n\ge1} n^{n-1}\frac{z^n}{n!}.
\end{equation}
Indeed, the (multivalued) function~$T(z)$ has singularities at
$0,\infty,e^{-1}$ only. 

\smallskip

A major theorem constrains the possible
growth of a holonomic function near any of its
singularities.
Paraphrasing Theorem~19.1 of~\cite[p.~111]{Wasow87}, we can
state:

\begin{theorem}[Structure Theorem for singularities]\label{classif-thm}
Let there be given a differential equation of the form~\eqref{h2}, a singular
point $z_0$, and a sector~$S$ with vertex at~$z_0$. Then,
for $z$ in a sufficiently narrow subsector~$S'$ of~$S$ and for
$|z-z_0|$ sufficiently small, there exists a basis of~$d$ linearly independent
solutions to~\eqref{h2}, such that
any solution~$Y$ in that basis admits as $z\to z_0$ 
in the subsector an asymptotic expansion of the form
\begin{equation}\label{struct-eqn}
Y\sim \exp\left(P(Z^{-1/r})\right)z^\alpha \sum_{j=0}^\infty Q_j(\log Z)Z^{js},
\qquad Z:=(z-z_0),
\end{equation}
where $P$ is a polynomial, $r$ an integer of~$\Z_{\ge0}$, $\alpha$ a
complex, $s$ a
rational of $\Q_{>0}$, and the $Q_j$ are a family of polynomials
of uniformly bounded degree. The
quantities $r,P,\alpha,s,Q_j$ depend on the particular solution and the
formal asymptotic expansions of~\eqref{struct-eqn} are $\C$-linearly
independent. 
\end{theorem}
(The argument is based on first constructing a \emph{formal} basis of
independent solutions, each of the form~(\ref{struct-eqn}), and then
applying to the possibly divergent expansions
a summation mechanism that converts such formal solutions into
actual analytic solutions. The restriction of the
statement to a subsector is related to
the Stokes phenomena associated to so-called
``irregular'' singularities.)

\smallskip

This theorem implies that the sequence $(n^{n-1}/n!)$ (hence $(n^n)$)
 is non-holonomic. Indeed the Cayley tree function satisfies the
functional equation
\[
T(z)=ze^{T(z)},
\]
corresponding to the fact that it enumerates labelled nonplane trees. 
Set $W(z)=-T(-z)$, which is otherwise known as the ``Lambert
W-function''. One has
\[
W(x)\mathop{=}_{x\to+\infty} \log x - \log\log x +O(1) ,\]
as verified by bootstrapping (see De Bruijn's
monograph~\cite[p.~26]{deBruijn81}). This is enough to conclude that
$W$, 
hence $T$, is non-holonomic as the $\log\log(\cdot)$ term is incompatible
with Eq.~(\ref{struct-eqn}). Observe that, conceptually, the 
proof involves considering the \emph{analytic continuation} of $T(z)$ 
and then extracting a clearly non-holonomic term in the expansion near a
singularity. More of this in the next sections.

Amongst other applications, we may cite:
\begin{itemize}
\item[---] Stanley's children rounds [\EIS{A066166}],
with exponential generating function
$(1-z)^{-z}$. The expansion as $z\to1$,
\[
(1-z)^{-z}
\mathop{\sim}_{z\to1} \frac{1}{1-z}\left(1+(1-z)\log(1-z)
+\frac{(1-z)^2\log^2(1-z)}{2!}+\cdots\right),
\]
contradicts the fact that logarithms can only appear with
bounded degrees in holonomic functions.
\item[---] Bell numbers have OGF $e^{e^z-1}$. In this case, the double
exponential behaviour as $z\to+\infty$ excludes them from the
holonomic ring.
\end{itemize}

\smallskip

Finally, what is given is often a sequence rather than a function.
Under such circumstances, it 
proves handy to be able to relate the asymptotic behaviour 
of~$f_n$ as $n\to+\infty$
to the asymptotic form of its generating function $f(z)$, near a
singularity.
Such transfers exists and are widely known in the literature as \emph{Abelian
theorems}.
We make use here of well-established principles in this theory,
as found, e.g., in the reference book by Bingham, Goldie, and
	Teugels~\cite{BiGoTe89}. 
For convenience of exposition, we state explicitly one version 
used repeatedly here:

\begin{theorem}[Basic Abelian theorem]\label{abel-thm}
Let $\phi(x)$ be any of the
functions
\begin{equation}\label{condabel}
x^\alpha (\log x)^{\beta}(\log\log x)^\gamma,\qquad\alpha\ge0,
\quad \beta,\gamma\in\C.
\end{equation}
Let $(u_n)$ be a sequence that satisfies the asymptotic estimate
\[
u_n \mathop{\sim}_{n\to\infty} \phi(n).
\]
Then, the generating function,
\[
u(z):=\sum_{n\ge 0} u_n z^n,
\]
satisfies the asymptotic estimate 
\begin{equation}\label{estim}
u(z)\mathop{\sim}_{z\to1-}
\Gamma(\alpha+1)\frac{1}{(1-z)}\phi\left(\frac{1}{1-z}\right).
\end{equation}
This estimate remains valid when $z$ tends to~1 in any sector with vertex
at~$1$, symmetric about the horizontal axis, and with opening angle
$<\pi$.
\end{theorem}
\begin{proof}[Proof (sketch)]
We shall content ourselves here with brief indications since
Corollary 1.7.3 p.~40 of~\cite{BiGoTe89}
provides simultaneously the needed Abelian
property and its real-analysis Tauberian converse\footnote{%
	The singularity analysis technology of Flajolet and
	Odlyzko~\cite{FlOd90b,Odlyzko95} provides sufficient
	conditions for the converse 
	complex-Tauberian implication.}, 
at least in the case when $z$
tends to~$1^-$ along the real axis. 

For simplicity, consider  first the representative
case where $\phi(x)=\log\log x$
and one has exactly $u_n=\phi(n)$ for $n\ge2$, with
$u_0=u_1=0$. Assume at this stage that~$z$ is
real positive and set $z=e^{-t}$, where $t\to0$ as $z\to1$.
We have
\[
u(z) = \sum_{n\ge2} \phi(n) e^{-nt}. \]

Take $n_1=\lfloor t^{-1}/\log t^{-1}\rfloor$. Basic majorizations imply
that the sum of
the terms corresponding to $n<n_1$ is bounded
from above by $n_1\log\log n_1$, which is
smaller than the right hand side of~(\ref{estim}).
Similarly, define $n_2=\lfloor t^{-1}\log t^{-1}\rfloor $.
The sum of terms with $n>n_2$ is easily checked to be $O(1)$. 
The remaining ``central'' terms $n_1\le n\le n_2$ are such that $\phi(n)$ varies
slowly over the interval and one has $\phi(n_1)\sim\phi(n_2)\sim
\phi(1/t)$.
One can thus take out a factor of $\phi(1/t)$
and conclude, upon approximating the sum by an integral, that
\begin{equation}\label{abelard}
\sum_{n=n_1}^{n_2} \phi(n) e^{-nt}\sim 
\frac{\phi(1/t)}{t}\int_{1/\log t^{-1}}^{\log t^{-1}} e^{-x}\, dx \sim 
\frac{\log\left(\log(1-z)^{-1}\right)}{1-z}.
\end{equation}
(Use the Euler-Maclaurin summation formula, then complete the tails.)

The proof above applies when 
\[
z=e^{-t+i\theta},\qquad
\hbox{with}\quad
|\theta|<\theta_0,\]
for some $\theta_0<\frac{\pi}{2}$. Once more only the central terms
matter asymptotically; the integral is then to be taken along a line
of angle $\theta$, but it reduces to the corresponding
 integral along the positive real line,
by virtue of the residue theorem. This suffices to justify the extension 
of the estimate to sectors. 
The case $\phi(x)=\log\log x$ is
then settled.

The extension to $u_n=\phi(n)$
when $\phi$ only involves powers of logarithms 
and of iterated logarithms follows similar lines, as $\log
n$ is also of slow variation. The inclusion of a power of
$n$, in the form $n^{\alpha}$ implies that the integral
in~(\ref{abelard}) should be modified to include a factor $t^{\alpha}$ 
leading to a real line integral that evaluates to the Gamma function,
$\Gamma(\alpha+1)$.

Finally, simple modifications of the previous arguments show that if a
sequence $(v_n)$ satisfies $v_n=o(\phi(n))$, for $\phi(n)$ any
of the sequences of~(\ref{condabel}), then its generating function~$v(z)$ is 
a little-oh of the right hand side of~(\ref{estim}). Decomposing $u_n=\phi(n)+v_n$
completes the proof of the statement.
\end{proof}

\begin{note} We have chosen to state the Abelian Theorem
(Theorem~\ref{abel-thm}) for $z$
varying in a cone of the complex plane, rather than the more customary
real line. In this way we can avail ourselves of the comparatively
simple Structure
Theorem (as stated above in Theorem~\ref{classif-thm})
and avoid some of the possible hardships due
to the Stokes phenomenon.
\end{note}

Here is a direct application of Theorem~\ref{abel-thm}. 
Let $\pi(x)$ be the number of primes less
than or equal to~$x$. By the Prime Number Theorem, one knows that
\[
\pi(n)\sim \frac{n}{\log n}.\]
The Abelian Theorem permits us to conclude about the
non-holonomic character of the sequence $(\pi(n))$, since 
\[
\sum_{n\ge1} \pi(n) z^n \mathop{\sim}_{z\to 1^-} \frac{1}{(1-z)^2
\log(1-z)^{-1}},\]
which contradicts what the Structure Theorem permits.

\begin{note}
In this article, we concentrate on proofs of non-holonomicity based on 
analysis, that is, eventually, \emph{asymptotic approximations}.
On a different register, powerful algebraic tools can be 
put to use in a number of situations.
Considerations on power series have been used 
by Harris and Sibuya~\cite{HaSi85} to show:
\emph{The reciprocal $(1/f)$ of a holonomic function~$f$ is holonomic
if and only  
if~$f'/f$ is algebraic}.
For instance, this proves the non-holonomicity of the reciprocal of 
Gauss' ${}_2F_1$ hypergeometric, except in degenerate cases.
Using differential Galois theory, Singer generalized this result 
in~\cite{Singer86}. He characterized the
possible polynomial relations between holonomic functions and also 
showed the following: \emph{A holonomic function~$f$ has
to be algebraic if any of $\exp\int f$ or $\phi(f)$ is holonomic, with 
$\phi$ an algebraic function of genus~$\ge1$}.
An analogous result for sequences is given in~\cite[Chap.~4]{SivdP97}:
\emph{If both~$f_n$ and~$1/f_n$ are holonomic, then $f_n$ is an interlacing 
of hypergeometric sequences}.
\end{note}

\section{The logarithmic sequence}

Let $f_n=\log n$ (with $\log 0\equiv0$) and let 
$f(z)$ be its generating function,
\[
f(z)=\sum_{n\ge 1} (\log n) z^n.\]
We propose to show that a sequence derived from $f_n$ by means of
holonomicity preserving transformations is non-holonomic. 
Consider a variant of the $n$th difference of the sequence $f_n$,
namely
\[
\widehat f_n:=\sum_{k=1}^n \binom{n}{k} (-1)^k \log k,\]
whose ordinary generating function $\widehat f(z):=\sum_{n\ge1} \widehat f_n
z^n$ has positive radius of convergence and satisfies
\begin{equation}\label{e2}
\widehat f(z) = \frac{1}{1-z}f\left(-\frac{z}{1-z}\right).
\end{equation}
It is known that holonomic functions are closed under product and
algebraic (hence also, rational) substitutions. Thus, $f$ and
$\widehat f$ are such that either both of them are holonomic
or none of them is holonomic. (See also
Stanley's paper~\cite[p.~181]{Stanley80} for a discussion of the 
fact that differencing preserves holonomicity.)

Next, Flajolet and Sedgewick proved in~\cite{FlSe95} that the sequence
$\widehat f_n$ satisfies the asymptotic estimate
\begin{equation}\label{e1}
\widehat f_n = \log\log n+O(1).
\end{equation}
As a matter of fact, a full expansion is derived
in~\cite[Th.~4]{FlSe95}, based on the N\"orlund-Rice integral
representation~\cite{Norlund54} 
\begin{equation}\label{norlund}
\widehat f_n = \frac{(-1)^n}{2i\pi}\int_{\mathcal H}
(\log s) \frac{n!}{s(s-1)\cdots (s-n)}\, ds,
\end{equation}
the use of a Hankel contour~$\mathcal H$, and estimates akin to those used 
for the determination of inverse Mellin transforms affected with an
algebraic-logarithmic singularity~\cite{Doetsch55}.

The proof of Proposition~\ref{log-thm} can now be easily completed.
By the Abelian estimate of Theorem~\ref{abel-thm} 
applied to the asymptotic form~(\ref{e1}) of $\widehat f_n$, we 
have
\[
\widehat f(z)\sim\frac{\log\left(\log (1-z)^{-1}\right)}{1-z}\qquad (z\to 1,~z\in S),\]
this in the whole of a sector $S$ 
with vertex at~1 and of
opening angle $<\pi$ extending towards the negative real axis 
symmetrically about the real axis.
Assume \emph{a contrario} that $f(z)$ is holonomic. Then $\widehat f(z)$, which
is associated to differences, is also holonomic. But then, given the 
Structure Theorem, a log-log asymptotic element 
valid in a subsector~$S'$
can never result from a $\C$-linear combination of elements,
each having the form~(\ref{struct-eqn}). A contradiction has thus been
reached, and Proposition~\ref{log-thm} is established.

\begin{note}
We have presented our proof in a way that seems to depend
on the imported estimate~(\ref{e1})
 of the logarithmic
differences. 
In this way, we could
save a few analytic steps. A conceptually equivalent and self-contained
proof would 
proceed from the asymptotic behaviour of
the analytic continuation of $f(z)$ as
$z\to-\infty$. (See our earlier discussion of~$T(z)$ for
a similar situation.) This can be achieved 
directly by means of a Lindel\"of integral
representation,
\[
f(-z)=\frac{1}{2i\pi}\int_{1/2-i\infty}^{1/2+i\infty} \left(\log
s\right) z^{s} \frac{\pi}
{\sin \pi s}\, ds.\]
(See Lindel\"of's monograph~\cite{Lindelof05} for explanations from the
mouth of the master and Flajolet's paper~\cite{Flajolet99} for related
developments.) It can then be verified,
by deforming the line of integration into a Hankel contour,
that non-holonomic elements crop up 
in the asymptotic expansion of~$f(z)$ at~$-\infty$.
\end{note}

\section{The sequence of powers}

The proof of  Proposition~\ref{pow-thm}
relies once more on the consideration of diagonal
differences. It has been established in~\cite{FlSe95} that 
\[
w_n:=\sum_{k=1}^n \binom{n}{k}(-1)^k k^\alpha,
\]
satisfies, for $\alpha\in\C\setminus\Z$,
\[
w_n = \frac{(\log n)^{-\alpha}}{\Gamma(1-\alpha)}
\left(1+O\left(\frac{1}{\log n}\right)\right).
\]
For instance,
\[
\sum_{k=1}^\infty \binom{n}{k}(-1)^k \sqrt{k}
=\frac{1}{\sqrt{\pi}{\log n}}+O\left((\log n)^{-3/2}\right).
\]
By the Abelian theorem, this implies for instance that the generating
function $w(z)$ of $w_n$ satisfies, for $\alpha=\frac12$ as $z\to1^{-}$ 
\[
w(z)\sim \frac{1}{\pi}\frac{1}{\sqrt{\log(1-z)^{-1}}}\frac{1}{1-z},
\]
while, for general $\alpha\not\in \Z$, the dominant asymptotics involve a
factor of the form 
\[
\left(\log(1-z)^{-1}\right)^{-\alpha}.\]
Such a factor is foreign to what the Structure Theorem provides as legal
holonomic asymptotics, as soon as $\alpha$ is nonintegral. This
completes the proof of Proposition~\ref{pow-thm}.

Proposition~\ref{pow-thm} implies in particular that the power sequences,
\[
n^{\sqrt{17}}, \quad n^i=\cos\log n+i \sin \log n,\quad
n^\pi,\ldots,
\]
are non-holonomic. 
Note that some of these sequences \emph{do} occur as
valid \emph{asymptotic} approximations of holonomic sequences, 
some of which even appear in natural combinatorial problems.
For instance, it is proved in~\cite{FlGoPuRo93} that 
the expected cost of a partial match
in a quadtree is holonomic and has the asymptotic form
$n^{(\sqrt{17}-3)/2}$. 

\section{The $n$th prime function}

Our proof of Proposition~\ref{prime-thm} will similarly involve detecting,
in the generating function $g(z)$ associated to primes, some elements 
that are incompatible
with holonomy and contradict the conclusions of
the Structure Theorem.

The $n$th prime function $n\mapsto g_n$ (often also written $p(n)$) is
a much researched function. In a way, this function is an inverse of
the function $\pi(x)$ that gives the number of primes in the interval
$[1,x]$. By the Prime Number Theorem, the function $\pi(x)$ satisfies
\[
\pi(x)=\Li(x)+R(x),\]
where $R(x)$ is of an order smaller than the main term
and $\Li(x)$ is the  logarithmic integral,
\[
\Li(x)=\int_2^x \frac{dt}{\log t} \sim \frac{x}{\log x}
\left(1+\frac{1!}{\log x}+\frac{2!}{(\log x)^2}+\cdots\right).\]
The precise 
description of the remainder term
$R(x)$ depends upon the Riemann hypothesis. However, it is known
\emph{unconditionally} that $R$ is small enough that the relation
$\pi(x)=y$ can be inverted asymptotically by just inverting the main
term~$\Li(x)$. In this way, one obtains (see~\cite{Salvy94} and
references therein) an estimate due to Cipolla~\cite{Cipolla02}, 
\begin{equation}\label{p1}
g_n = n\log n +n\log\log n +O(n).
\end{equation}
The log-log term is once more a barrier to holonomicity.

To see this, we introduce now the function,
\[
\ell(z)=\frac{z}{(1-z)^2}\log\frac{1}{1-z}+\frac{z}{(1-z)^2},
\]
which is clearly holonomic.
It satisfies, with $H_n=1+\frac12+\cdots\frac{1}{n}$ denoting the
harmonic number, 
\[
\ell_n=nH_n = n\log n +O(n).
\]
Then, by taking a difference, 
 one gets
\[
g_n-\ell_n = n\log\log n +O(n).\]
Thus, the difference $g(z)-\ell(z)$ satisfies
\[
g(z)-h(z) \mathop{\sim}_{z\to1^-}
\frac{\log\left(\log(1-z)^{-1}\right)}{(1-z)^2}.
\]
This last fact is incompatible with holonomicity, by the very same argument
as in the previous section.

\section{Conclusion}

As we have strived to illustrate and as is
otherwise seen in many areas of mathematics (see the opening quotation),
singularities are central to 
the understanding of properties of numeric sequences. 
It should be clear by now that a large number of sequences can be
proved to be non-holonomic. Here is a brief recapitulation of methods 
and possible extensions.

\medskip
{\bf1. \emph{Asymptotic discrepancies}}.
In order to conclude
that a sequence $(u_n)$ is non-holonomic,
the following conditions are sufficient: \emph{$(i)$~the generating function
of the sequence
$(u_n)$ (or one of its cognates)
admits, near a singularity, an asymptotic expansion in a scale that
involves logarithms and 
iterated logarithms; $(ii)$~at least one term in that expansion 
is an iterated
logarithm or  a 
power of a logarithm with an exponent not in~$\Z_{\ge0}$.}
It is then apparent from our earlier developments that sequences like
\[
\sqrt{n^7+1},\qquad \frac{1}{H_n},\quad
\sqrt{\frac{\log n}{H_n}},\quad \log \frac{p(2n)}{p(n)},\quad
\frac{n}{\sqrt{n}+\log n},\quad 
p(n^2),\]
($H_n$ the harmonic number and $p(n)$ the $n$th prime function)
fail to be holonomic.
Also, techniques of this note extend easily to other slowly varying
sequences, like 
\[
e^{\sqrt{\log n}}, \quad 
\log\log\log n, \ldots\,,\]
as well as to any sequence that involves any such term somewhere in its
asymptotic expansion.

The singularity-based technology has otherwise been used to establish the
non-algebraic character of sequences arising from combinatorics and
the theory of formal languages in~\cite{Flajolet87b}. Once more, such 
transcendence results imply that, structurally, the corresponding
objects cannot be (unambiguously) encoded by words of a context-free
language. For instance, two-dimensional walks on  a regular
lattice that are constrained to the first quadrant cannot be
described (via a length-preserving encoding) by means of an  
unambiguous context-free grammar. This property is neatly visible from a
logarithmic component in the generating function of walks~\cite{FlSe04}, which
contradicts the Structure Theorem for algebraic functions (also known
as Newton-Puiseux!).

\medskip
{\bf2. \emph{Infinitude of singularities}}.
A famous theorem of P\'olya and Carlson~\cite{Bieberbach31} implies
the following\footnote{%
	Stefan Gerhold is grateful to Richard Stanley for pointing out
	the P\'olya-Carlson connection. 
}: 
\emph{A function analytic at the origin having integer
coefficients and assumed to converge in the open unit disc is either a
rational function or else it admits the unit circle as a natural
boundary.}  Consider then the generating function $g(z)$ of $g_n\equiv
p(n)$, the $n$th prime function, which can be subjected to
P\'olya-Carlson. Either it has a natural boundary, in which case it
cannot be holonomic, since holonomic functions have isolated
singularities. Else, it is rational; but this would be a clear
contradiction, since no rational function can have coefficients of the
asymptotic form $n\log n$ (by virtue of
the Abelian Theorem, say). This proof\footnote{%
	Alternatively, Erd\H os, Maxsein, and Smith~\cite{ErMaSm90}
	have shown that an integer 
	recurrent sequence [i.e., one with a rational generating function] 
	which consists only of primes is 
	necessarily a periodic sequence, therefore involving only finitely
	many different values. We are indebted to Dr Alin Bostan for
	this remark.} 
seems to require the Prime Number Theorem but many weaker bounds that
are elementary (e.g., the ones due to Chebyshev) are sufficient for
this purpose. It is  pleasant to note that the P\'olya-Carlson Theorem
was earlier employed in a similar fashion~\cite{AuFlGa87} in order
to establish a structural lower bound in the theory of formal languages.

\medskip
{\bf 3. \emph{Arithmetic discrepancies\footnote{%
	See the recent monograph of Everest \emph{al.}~\cite{EvPoSh03}
	for a compendium of results relative to recurrence sequences.
}.}}
In~\cite[p.~294]{Flajolet87b}, a sketch was
given of the solution to a conjecture of Stanley~\cite{Stanley80} 
to the effect that the generating function 
\[
S_k(z)=\sum_{n\ge0} \binom{2n}{n}^k z^n,
\]
is transcendental for odd
values of $k$, $k=3,5,\ldots$. 
(As already noted by Stanley,
 	$S_k$ is clearly transcendental
	for~$k$ even because of the presence of logarithmic factors.)
Longer algebraic proofs have since  been published,
see~\cite{WoSh89} and~\cite{Allouche97} for a discussion.
 It
suffices to remark, as shown by an Abelian argument, that the local expansion
at the finite singularity of~$S_k$ for $k=2\ell+1$ odd, satisfies
\begin{equation}\label{stan1}
\frac{d^\ell}{dz^\ell}
S_{2\ell+1}\left(\frac{z}{4^{2\ell+1}}\right) \mathop{\sim}_{z\to1^{-}}
c_\ell\pi^{-\ell}\frac{1}{\sqrt{1-z}}\,\qquad c_\ell\in \Q.
\end{equation}
This asymptotic form is incompatible with algebraicity.
Indeed, if $S_k$ were algebraic over $\C(z)$, it would be algebraic
over $\Q(z)$, as it has rational Taylor coefficients (this, by a famous lemma
of Eisenstein~\cite{Bieberbach31}). But in that case, its Puiseux expansion
could only involve  algebraic numbers. Equation~(\ref{stan1})
contradicts this, by virtue of the transcendence of~$\pi$. 
Et voila!
In this last case, we have to play not only with the \emph{shape} of
an asymptotic expansion, but also with the \emph{arithmetic nature} of its
coefficients.

\def\cprime{$'$}
\providecommand{\bysame}{\leavevmode\hbox to3em{\hrulefill}\thinspace}
\providecommand{\MR}{\relax\ifhmode\unskip\space\fi MR }
\providecommand{\MRhref}[2]{%
  \href{http://www.ams.org/mathscinet-getitem?mr=#1}{#2}
}
\providecommand{\href}[2]{#2}

\end{document}